\documentclass[MNA]{mna} % 

\usepackage[T1]{fontenc}
\usepackage[latin1]{inputenc}

\usepackage[T1]{fontenc}
\usepackage[latin1]{inputenc}
\usepackage[english]{babel}
\SHORTTITLE{Two-layer neural model}

\TITLE{Mean field system of a two-layer neural model in a diffusive regime}
\AUTHORS{%
  Xavier~Erny\footnote{Universit\'e Paris-Saclay, CNRS, Univ Evry, Laboratoire de Math\'ematiques et Mod\'elisation d'Evry, 91037, Evry, France
    \BEMAIL{xav.erny@orange.fr}}}

\SHORTAUTHOR{Xavier Erny}

\AMSSUBJ{60J35 ; 60J60 ; 60K35 ; 60G55} % Edit. Separate items with ;
\SUBMITTED{June 10, 2021} % Edit.
\ACCEPTED{April 13, 2022} % Edit.

\ARXIVID{2106.05008} % Edit.
\HALID{hal-03255276} % Edit.

\VOLUME{2}
\YEAR{2022}
\PAPERNUM{3}
\DOI{10.46298/mna.7570}

\ABSTRACT{We study a model of interacting neurons. This neural system is composed of two layers of neurons such that the neurons of the first layer send their spikes to the neurons of the second one: if $N$ is the number of neurons of the first layer, at each spiking time of the first layer, every neuron of both layers receives an amount of potential of the form $U/\sqrt{N},$ where $U$ is a centered random variable. This kind of structure of neurons can model a part of the structure of the visual cortex: the first layer represents the primary visual cortex V1 and the second one the visual area V2. The model consists of two stochastic processes, one modelling the membrane potential of the neurons of the first layer, and the other the membrane potential of the neurons of the second one. We prove the convergence of these processes as the number of neurons~$N$ goes to infinity and obtain a convergence speed. The proofs rely on similar arguments as those used in \cite{erny_mean_2022}: the convergence speed of the semigroups of the processes is obtained from the convergence speed of their infinitesimal generators using a Trotter-Kato formula, and from the regularity of the limit semigroup.}

\newcommand\n{\mathbb{N}}
\renewcommand\r{\mathbb{R}}

\renewcommand\ll{\left}
\newcommand\rr{\right}
\newcommand\un{\mathds{1}}

\renewcommand\phi{\varphi}
\newcommand\pro[1]{\mathbb{P}\left(#1\right)}
\newcommand\esp[1]{\mathbb{E}\left[#1\right]}

\newcommand\espcc[2]{\mathbb{E}_{#1}\left[#2\right]}

\newcommand\uno[1]{\un_{\left\{#1\right\}}}
\newcommand\1{\un}

\begin{document}

%%%%%%%%%%%%%%%%%%%%%%%%%%%%%%%%%%%%%%%%%%%%%%%%%%%%%%%%%%%%%%%%%%%
%%                                                               %%
%% No need for \maketitle.                                       %%
%%                                                               %%
%%%%%%%%%%%%%%%%%%%%%%%%%%%%%%%%%%%%%%%%%%%%%%%%%%%%%%%%%%%%%%%%%%%

%%%%%%%%%%%%%%%%%%%%%%%%%%%%%%%%%%%%%%%%%%%%%%%%%%%%%%%%%%%%%%%%%%%
%%                                                               %%
%% Please replace what follows by the body of your article       %%
%% (up to the bibliography):                                     %%
%%                                                               %%
%%%%%%%%%%%%%%%%%%%%%%%%%%%%%%%%%%%%%%%%%%%%%%%%%%%%%%%%%%%%%%%%%%%

\section{Introduction}

\subsection{Model}

The model we study in this paper is a mean-field model of interacting neurons. Such models have been studied widely in the literature to model the dynamics of neural networks (see e.g. \cite{faugeras_constructive_2009, baladron_mean-field_2012,lucon_mean_2014, de_masi_hydrodynamic_2015} and the references therein). A natural way of modelling the activity of the neurons is to use point processes to represent their spike trains. In this context, Hawkes processes have called attention as a natural model for the spike trains of neurons. These processes have been introduced in \cite{hawkes_spectra_1971} to model the aftershocks of earthquakes. By definition, Hawkes processes are mutually exciting point processes, therefore they can be naturally used to model the mutually exciting/inhibiting behavior of neurons (see \cite{pillow_spatio-temporal_2008, grun_analysis_2010, reynaud-bouret_goodness--fit_2014}). Some authors have also managed to model the repolarization of the neurons using Hawkes processes with variable length memory (see \cite{fournier_toy_2016, hodara_hawkes_2017, erny_conditional_2021}). Finally "age dependent Hawkes processes" introduced in \cite{chevallier_mean-field_2017} and further studied in \cite{chevallier_fluctuations_2017, bonde_raad_stability_2020} can be used to make the spike rate of a neuron  depend on the delay since its last spike, allowing to consider for example a refractory period for the neurons.

This paper aims to study a model of two layers of neurons such that:
\begin{itemize}
	\item the neurons of the first layer send their spikes to the neurons of both the first and second layers,
	\item the spikes of the neurons of the second layer do not affect these two layers.
\end{itemize}

We note $N$ the number of neurons of the first layer. The idea of our modelling is to study the "averaged effect" of the first layer on a single neuron of the second one. In the following, we model this averaged effect through the averaged value of the membrane potentials of the neurons of the first layer. Let us note $X^N_t$ the averaged value of these potential at time~$t$, and $Y^N_t$ the membrane potential of one particular neuron of the second layer at time~$t$.

We assume that the dynamics of the processes $X^N$ and $Y^N$ are given by the following SDEs:
\begin{align}
X^N_t=&X^N_0-\alpha_1 \int_0^tX^N_sds + \frac{1}{\sqrt{N}}\sum_{j=1}^{N}\int_{[0,t]\times\r_+\times\r}u\cdot \uno{z\leq f_1(X^N_{s-})}d\pi^{(1),j}(s,z,u)\label{XNYN}\\
&-\frac1N\sum_{j=1}^{N}\int_{[0,t]\times\r_+\times\r}X^{N}_{s-}\uno{z\leq f_1(X^N_{s-})}d\pi^{(1),j}(s,z,u),\nonumber\\
Y^N_t=&Y^N_0-\alpha_2 \int_0^tY^N_sds + \frac{1}{\sqrt{N}}\sum_{j=1}^{N}\int_{[0,t]\times\r_+\times\r}u\cdot \uno{z\leq f_1(X^N_{s-})}d\pi^{(1),j}(s,z,u)\nonumber\\
& - \int_{[0,t]\times \r_+}Y^N_{s-}\uno{z\leq f_2(Y^N_{s-})}d\pi^{(2)}(s,z),\nonumber
\end{align}
where $X^N_{t-}$ (resp. $Y^N_{t-}$) denotes the limit from the left of $X^N$ (resp. $Y^N$) at time~$t.$ Above, $\alpha_1,\alpha_2$ are positive parameters (resp. $f_1,f_2$ are non-negative functions) of the model that represent the leakage rate characteristics (resp. the spike rate functions) of respectively the neurons of the first layer and the neurons of the second one. In~\eqref{XNYN}, $\pi^{(1),j}$ are i.i.d. Poisson measures with intensity $ds\cdot dz\cdot d\nu(u)$ ($\nu$ is a centered probability measure on~$\r$), and $\pi^{(2)}$ a Poisson measure with Lebesgue intensity independent of the other Poisson measures.

	\begin{remark}
		Note that it is possible to rewrite equation~\eqref{XNYN} in a more compact form:
		\begin{align*}
		X^N_t=&X^N_0-\alpha_1 \int_0^tX^N_sds + \frac{1}{\sqrt{N}}\int_{[0,t]\times\r_+\times\r}u\cdot \uno{z\leq f_1(X^N_{s-})}d\pi^{(1)}(s,z,u)\\
		&-\frac1N\int_{[0,t]\times\r_+\times\r}X^{N}_{s-}\uno{z\leq f_1(X^N_{s-})}d\pi^{(1)}(s,z,u),\\
		Y^N_t=&Y^N_0-\alpha_2 \int_0^tY^N_sds + \frac{1}{\sqrt{N}}\int_{[0,t]\times\r_+\times\r}u\cdot \uno{z\leq f_1(X^N_{s-})}d\pi^{(1)}(s,z,u)\\
		& - \int_{[0,t]\times \r_+}Y^N_{s-}\uno{z\leq f_2(Y^N_{s-})}d\pi^{(2)}(s,z),
		\end{align*}
		where $\pi^{(1)}$ is a Poisson measure with intensity $N\cdot ds\cdot dz\cdot d\nu(u).$ It is formally possible to define $\pi^{(1)}$ as the union of the $\pi^{(1),j}$ ($1\leq j\leq N$). In the rest of the paper, we use the form~\eqref{XNYN} because it is more understandable for the model.
	\end{remark}

In~\eqref{XNYN}, the counting processes
$$P^{(1),N,i}_t := \int_{[0,t]\times\r_+\times\r}\uno{z\leq f_1(X^N_{s-})}d\pi^{(1),i}(s,z,u),~~1\leq i\leq N,$$
and
$$P^{(2),N}_t := \int_{[0,t]\times\r_+}\uno{z\leq f_2(X^N_{s-})}d\pi^{(2)}(s,z),$$
model respectively the spike trains of the neuron~$i$ of the first layer ($1\leq i\leq N$) and that of the particular neuron of the second one. In this model, the form of the synaptic weights is $u/\sqrt{N}$ where $u$ is the third variable of the Poisson measures~$\pi^{(1),i}$ ($1\leq i\leq N$). These variables can be seen as random variables that are i.i.d. $\nu-$distributed, and whence centered. Note that it is necessary to consider centered synaptic weights because of the normalization~$N^{-1/2},$ and in order to make appear these centered variables~$u$, we need to write the point processes $P^{(1),N,i}$ as thinnings of Poisson measures in~\eqref{XNYN}.

Let us explain in other words the dynamics of the equations~\eqref{XNYN}. The second term of both SDEs describes the effect of the spikes of the neurons of the first layer on both layers. The last term in the SDE of $Y^N$ (i.e. the fourth line of~\eqref{XNYN}) models the repolarization of the potential of the particular neuron of the second layer: right after the neuron emits a spike, its potential goes fast (here it jumps immediately) to its resting value that we assume to be zero. The last term of the SDE of $X^N$ (i.e. the second line of~\eqref{XNYN}) should be interpreted as an averaged repolarization. Indeed, the neurons of the first layers fire at rate $N\cdot f_1(X^N_{t-}).$ So $X^N_t$ should be interpreted as the averaged potential of the neurons of the first layer. And, according to our model, when a neuron of the first layer sends a spike, its potential should jump to zero. As a consequence, if we assume that this potential was "close" to the averaged value $X^N_{t-}$ before the jump, then, right after the spike, the quantity lost by the averaged potential is $N^{-1}X^N_{t-}.$

The reason why we only consider the effect of the first layer of neurons through their averaged potential~$X^N$ is to guarantee our model~$(X^N,Y^N)$ to be two-dimensional, and not $(N+1)-$dimensional. In this second frame, the technics used in the proofs would fail. This if for the same reason that we assume our neural networks to be complete graphs, since it implies that the dynamics of all the neurons are similar.

The main result of this paper (Theorem~\ref{mainresult}) is the convergence in distribution as $N$ goes to infinity (with an explicit convergence speed) of the sequence of the two-dimensional Markov processes $(X^N,Y^N)_N$ (defined in~\eqref{XNYN}) to the following limit process
\begin{align}
\bar X_t=& \bar X_0 -\alpha_1\int_0^t\bar X_sds + \sigma\int_0^t\sqrt{f_1(\bar X_s)}dW_s - \int_0^t\bar X_s f_1(\bar X_s)ds,\label{barXbarY}\\
\bar Y_t=&\bar Y_0 - \alpha_2\int_0^t\bar Y_sds + \sigma\int_0^t\sqrt{f_1(\bar X_s)}dW_s - \int_{[0,t]\times\r_+} \bar Y_{s-} \uno{z\leq f_2(\bar Y_{s-})}d\bar\pi(s,z),\nonumber
\end{align}
where $W$ is a standard one-dimensional Brownian motion, $\bar\pi$ is a Poisson measure on $\r_+^2$ with Lebesgue intensity independent of~$W$, and $\sigma^2$ is the variance of the probability measure~$\nu$.

\begin{remark}\label{simulationinteret}
	Note that Theorem~\ref{mainresult} has a practical interest if one wants to simulate the trajectories of the process~$(X^N,Y^N).$ According to this result, the limit process $(\bar X,\bar Y)$ is an approximation of $(X^N,Y^N),$ and it is easier to simulate. Indeed, to simulate $(X^N,Y^N)$ one has to compute the trajectories of the $N$ Poisson measures $\pi^{(1),i}$ ($1\leq i\leq N$), which can be computationally hard if $N$ is large. So it can be easier to simulate a single Brownian motion~$W$ instead of those Poisson measures.
	
	Then, if we simulate the limit process $(\bar X,\bar Y)$ instead of the "real" one $(X^N,Y^N)$ (i.e. the one that actually corresponds to the neural model), we make an error. The interesting point of Theorem~\ref{mainresult} is to quantify this error (i.e. the convergence speed). Let us mention that this error is of order $N^{-1/2},$ which can be considered "small" if the first layer represents for example the primary visual cortex V1, since the number of neurons~$N$ in~V1 is estimated to be of order~$10^8$ according to \cite{leuba_changes_1994}.
\end{remark}

\subsection{Biological motivation}

This multi-layers structure is a natural model for the visual cortex, which is composed of five interacting and functionally distinct layers of neurons (see Chapter~5 of \cite{hubel_eye_1995} for a detailed description of the structure of the visual cortex). The two-layer model of this paper can be seen as a model of the primary visual cortex V1 (the first layer) and the visual area V2 (the second layer). A more interesting model could consist in considering five layers to represent~V1-V5 and more complex interactions between these layers. As we have mentioned earlier, to simplify the model we only consider an averaged potential $X^N$ for the neurons of the first layer. We also assume that all the neurons of the second layer have the same dynamics, and that is why we just study one particular neuron of the second layer.

In our model, we assume our first layer to be a complete graph (i.e. each neuron is connected by a synapse to any other neuron). As it is explained above, this property is important in our proofs, but it is a drawback of the biological model since a neuron has only a limited number of synaptic connections compared to the number of neurons of the network. It is possible to extend this model to the case where the first layer is divided into a fix number~$k$ of sublayers. In this case, each sublayer has to be assumed to be a complete graph and one can consider any interaction structure between the sublayers. In this new model, assuming that the number~$k$ of sublayers does not depend on~$N$, the results of the paper would hold true, and the techniques used to prove the results would be the same (the main difference is that $X^N$ would be a $k-$dimensional process instead of an one-dimensional one). This model may be interesting to model the six functionally distinct layers of V1, and even the four sublayers of the fourth one of V1 (i.e. taking $k=9$) (see e.g. Chapter~5 of \cite{hubel_eye_1995}). In order to ease the readability of the proofs, we prefer not to study the "multi-sublayers model". Let us also mention that for similar reasons, it is possible to extend our model to study a fix number of neurons of the second layer and not only one (provided that this number does not depend on~$N$).

Let us end this section by commenting the form of the synaptic weights $u/\sqrt{N}$ where the $u$ are random variables. These centered synaptic weights model the fact that we consider "balanced networks" (i.e. neural networks where the excitatory inputs and the inhibitory ones are balanced). Experiments support this hypothesis of balanced synaptic connections as a natural model (see \cite{shu_turning_2003, haider_neocortical_2006} and the references therein). Let us also mention that, the experiments done in \cite{barral_synaptic_2016} support the assumptions that the order of magnitude of the synaptic weights is close to the theoretical value (i.e. the square root of the number of presynaptic neurons, since the interactions are balanced), which is $\sqrt{N}$ with our assumptions.

\subsection{Mathematical techniques}

Let us give a brief overview of the mathematical studies of balanced network: \cite{faugeras_asymptotic_2015} have proved a large deviation principle for a discrete-time model with correlated synaptic weights, in \cite{sompolinsky_chaos_1988} the authors have studied a continuous-time model of ODEs parametrized by random synaptic weights where at the limit a Gaussian field is created by these weights, \cite{rosenbaum_balanced_2014} have studied the exsitence of a balanced state for the large scale limit of a similar model, and \cite{pfaffelhuber_mean-field_2021} have proved the convergence of Hawkes processes in distribution to a limit process that is solution to a stochastic convolution equation (also known as stochastic Volterra equation of convolution type). The form of the model of \cite{pfaffelhuber_mean-field_2021} is close to the one we consider here. In these papers, the synaptic weight of a given synapse is a random variable and each time a neuron sends a spike to another neuron, the potential of this other neuron always receives the same random variable from this synapse. In our model (similarly as in \cite{erny_mean_2022, erny_conditional_2021}), at each spiking time, we make a choice of random synaptic weight independently of the previous ones. This implies that the role of a synapse (i.e. excitatory or inhibitory) can change.
This property is inconsistent with biological evidence, but compared to the model of \cite{pfaffelhuber_mean-field_2021}, it allows us to obtain a convergence speed in addition to the convergence in distribution of the model, and it also allows us to consider more general distributions than just Rademacher distribution for the law of the synaptic weights.

Let us also note that this kind of multi-layers model has already been studied in \cite{ditlevsen_multi-class_2017} and at Section~5 of \cite{erny_white-noise_2021}. In both these references, the number of layers is not limited to two (but does not depend on~$N$) and every membrane potential is modeled by a specific process without the "repolarization" effect in the equations. Consequently the dimension of the $N-$neurons model is~$N$, and as we have claimed previously, the techniques of this paper would not work on these models. In \cite{ditlevsen_multi-class_2017}, the authors have obtained a $L^1-$convergence with a convergence speed in $N^{-1/2}.$ The processes they have studied are Hawkes processes and the order of synaptic weight is $N^{-1}.$ The model in \cite{erny_white-noise_2021} is closer to~\eqref{XNYN} since the synaptic weight is also centered and of order~$N^{-1/2}.$ The advantage of the model of the current paper is that we obtain a convergence speed of order~$N^{-1/2}$ whereas in \cite{erny_white-noise_2021}, we have just proved a convergence in distribution and in \cite{erny_strong_2021} (where the results hold true in a multi-layers framework) the order of the convergence speed is~$N^{-1/10}.$ Note that it is interesting to have a good convergence speed for simulation purpose (see Remark~\ref{simulationinteret}).

The question of the convergence of system with this kind of repolarization effect has already been studied previously (in a one-layer framework): see \cite{de_masi_hydrodynamic_2015, andreis_mckeanvlasov_2018, erny_well-posedness_2021} for models with synaptic weights of order~$N^{-1}$, and \cite{erny_conditional_2021} for synaptic weights of order~$N^{-1/2}$. The novelty of this article consists in obtaining a convergence speed for the convergence in distribution where the order of the synaptic weights is~$N^{-1/2}.$

We can also remark that~\eqref{XNYN} is a generalization of the model of \cite{erny_mean_2022} where we have studied the convergence of the one-dimensional process~$X^N$ in a simplified framework (i.e. without the last term in the SDE of $X^N$). The interesting points of the model of this paper is that it is more relevant for the modelling of a particular structure of neural network (e.g. the layers V1 and V2 of the virtual cortex), and that we can still prove similar results.

In order to quantify the convergence speed of the distribution of $(X^N,Y^N)$ to the distribution of $(\bar X,\bar Y),$ we rely on their semigroups, respectively denoted by $(P^N_t)_t$ and $(\bar P_t)_t.$ By definition the semigroups are families of operators that can be defined on the space of continuous and bounded functions by the following: for $g$ continuous and bounded, and $x,y\in\r,$
$$P^N_tg(x,y) := \espcc{(x,y)}{g(X^N_t,Y^N_t)},$$
where $\mathbb{E}_{(x,y)}$ is the expectation related to the probability measure $\mathbb{P}_{(x,y)}$ under which the initial value of the process $(X^N,Y^N)$ is given by the deterministic value $(x,y)=(X^N_0,Y^N_0).$

To obtain this convergence speed, we use similar arguments as in \cite{erny_mean_2022}. We consider the infinitesimal generators $A^N$ and $\bar A$ of respectively $P^N$ and $\bar P.$ These infinitesimal generators are operators that can be defined as the derivatives of their respective semigroups:
$$A^Ng(x,y) = \underset{t\rightarrow 0}{\lim}\frac{1}{t}\ll(P^N_tg(x,y) - g(x,y)\rr).$$
The notion of generators can slightly differ from a reference to another for technical reasons. This is why we define precisely the notion that we use in this paper at the beginning of Section~\ref{generator} (it is the same notion that we have introduced and studied in Appendix~$A$ of \cite{erny_mean_2022}). 

We begin by proving the convergence of these generators and by obtaining a convergence speed. This convergence speed depends on the derivatives of the test-function~$g$. Then, we deduce a convergence speed for the semigroups using the following formula:
\begin{equation*}
\bar P_tg(x,y) - P^N_tg(x,y) = \int_0^t P^N_{t-s}\ll(\bar A - A^N\rr)\bar P_sg(x,y)ds.
\end{equation*}

This kind of techniques used to obtain the convergence of the semigroups from the convergence of the generators are similar to the results of the Section~1.6 of \cite{ethier_markov_2005}. In particular, the formula above is proved in Lemma~1.6.2 of \cite{ethier_markov_2005} with stronger assumptions on the generators compared to the result of this paper and~\cite{erny_mean_2022}.

To be able to use~\eqref{trotterkato} to deduce a convergence speed for the semigroups from the one for the generators, we need to control the differentiability of the function $(x,y)\mapsto\bar P_sg(x,y).$ So, one of the main step of the proof of our main result consists in studying the regularity of the limit semigroup $(\bar P_t)_t.$ For this purpose, we need to study the regularity of the related stochastic flow (see Section~\ref{notation} for the notion of stochastic flow). Compared to the proofs of \cite{erny_mean_2022}, there are two additional technical difficulties that we handle at Section~\ref{semigroups}:
\begin{itemize}
	\item The trajectories of the process $(\bar X,\bar Y)$ are not continuous. Whence, the regularity of its semigroup is harder to study.
	\item The jump times of the process $(\bar X,\bar Y)$ occur at the rate $f_2(\bar Y_{t-}).$ This implies that this rate depends on the initial condition of the process. The regularity of this dependency is hard to study directly.
\end{itemize}

To overcome the first point, we study the regularity of the stochastic flow between the jump times (see Proposition~\ref{flotregulier}). Let us remark that \cite{bally_regularity_2018} have already used the idea of working between the jump times of the flow to study its regularity and that of its semigroup. For the second point, we introduce an auxiliary limit process $(\tilde X,\tilde Y)$ for which the rate of the jumps is constant. Then, we deduce the regularity of the semigroup $(\bar P_t)_t$ from the regularity of the flow of $(\tilde X,\tilde Y)$ using Girsanov's theorem for jump processes that gives the explicit Radon-Nikodym derivative between the processes $(\bar X,\bar Y)$ and $(\tilde X,\tilde Y)$ (see Proposition~\ref{girsanov}).

\subsection{Notation}\label{notation}

Throughout the paper, we use the following notation:
\begin{itemize}
	\item For $n,d\in\n^*,$ we note $C^n_b(\r^d)$ the set of $\r-$valued functions defined on~$\r^d$ that are $C^n,$ bounded and such that all their partial derivatives up to order~$n$ are bounded.
		\item We call $\beta\in\n^2$ a multi-index and we note $|\beta| = \beta_1 + \beta_2$. If $g:\r^2\rightarrow\r$ is $C^n$ and if $\beta$ is a multi-index such that $|\beta|\leq n,$ we write $\partial_\beta g$ the partial derivative of $g$ w.r.t.~$\beta.$ For the first order partial derivatives, we write instead $\partial_1g$ and $\partial_2g,$ and for the second order derivatives, we write $\partial^2_{i,j}g$ ($i,j\in\{1,2\}$). When there is no possible confusion, we may use the classical notation $\partial_xg(x,y) = \partial_1g(x,y).$
	\item For $g\in C^n_b(\r^2),$ we note
	$$||g||_{n,\infty} = \sum_{|\beta|\leq n} ||\partial_\beta g||_\infty.$$
	\item Given a process $X$ (resp. a multi-dimensional process $(X,Y)$) that is solution to some SDE, we note $X^{(x)}$ (resp. $(X^{(x,y)},Y^{(x,y)})$) the process solution to the SDE w.r.t. the initial condition $X_0=x$ (resp. $(X_0,Y_0)=(x,y)$). The function $x\mapsto X^{(x)}$ (resp. $(x,y)\mapsto (X^{(x,y)},Y^{(x,y)})$) is called the stochastic flow of the process $X$ (resp. $(X,Y)$). By definition, a stochastic flow of a Markov process represents its trajectories as functions of its initial condition (see Definition~1.1.1 of \cite{kunita_lectures_1986} for a more general definition). In this paper, we speak about stochastic flow instead of the trajectories of the process itself when we want to emphasize the fact that it is a function of the initial value, and to use the same terminology as \cite{kunita_lectures_1986}.
	\item In all the paper, $C$ denotes an arbitrary positive constant that is independent of any parameters except the parameters of the model. If a constant depends on some parameter~$\theta,$ we write $C_\theta$ instead. The value of such constants can change from line to line in an equation.
\end{itemize}

\subsection{Assumptions and main result}

To begin with, let us state the conditions we need to study the model.
\begin{assumptions}\label{hyp}$ $
	\begin{itemize}
		\item[(i)] The function $f_1$ is $C^1_b(\r)$, positive, and, for some constant $C>0,$ for any $x\geq 0,$ $|f_1'(x)|\leq C/(1+|x|)$,
		\item[(ii)] the function $f_2$ is non-negative and belongs to~$C^3_b(\r)$,
		\item[(iii)] the function $\sqrt{f_1}$ belongs to~$C^4_b(\r)$,
		\item[(iv)] the probability measure $\nu$ is centered and admits a finite third moment. We note $\sigma^2$ its variance.
	\end{itemize}
\end{assumptions}

\begin{remark}
	Any sigmoid function satisfies the conditions imposed by Assump.~\ref{hyp}$.(i)-(iii)$ on the functions $f_1$ and $f_2.$ This kind of function has been shown to be natural for neural models (see \cite{velichko_oscillator_2020}).
\end{remark}

Under Assumption~\ref{hyp}$.(i)-(ii)$, the strong well-posedness of the SDEs~\eqref{XNYN} and~\eqref{barXbarY} is a classical result since these two conditions guarantee that the coefficients are Lipschitz continuous and the jump rates bounded (see Theorem~IV.9.1 of~\cite{ikeda_stochastic_1989}).

The boundedness of $f_1$ and the assumption on $f_1'$ are used to prove the well-posedness of the limit process~\eqref{barXbarY} and to simplify some technical details, since it implies that the function $x\mapsto xf_1(x)$ is Lipschitz continuous. On the other hand, the boundedness of $f_2$ is important to control the Radon-Nikodym derivative of the limit process~$(\bar X,\bar Y)$ with respect to the auxiliary process~$(\tilde X,\tilde Y)$ that we introduce at Section~\ref{semigroups} (see Proposition~\ref{girsanov} for the form of this derivative).

Condition~$(iii)$ allows to prove that the stochastic flow of the auxiliary process is $C^3$ between its jump times (see Proposition~\ref{flotc3}).
And condition~$(iv)$ allows us to obtain a convergence speed for the generators (see Lemma~\ref{convergencegenerator}).

Under the additional assumption that $(X^N_0,Y^N_0)$ converges in distribution to $(\bar X_0,\bar Y_0)$, it is known that the process $(X^N,Y^N)$ converges in distribution to $(\bar X,\bar Y)$ in Skorokhod topology (see Theorem~IX.4.21 of \cite{jacod_limit_2003}). The main result of this paper consists in obtaining the convergence speed of the semigroups $P^N$ to the semigroup $\bar P$ of these processes (recalling that the processes $(X^N,Y^N)$ and $(\bar X,\bar Y)$ are defined respectively at~\eqref{XNYN} and~\eqref{barXbarY}, and that a semigroup can be defined as $P^N_tg(x,y) := \mathbb{E}[g(X^N_t,Y^N_t)|(X^N_0,Y^N_0)=(x,y)]$).

\begin{theorem}
	\label{mainresult}
	Grant Assumptions~\ref{hyp}$.(i)-(iv)$. For all $t\geq 0,$ $N\in\n^*,$ $g\in C^3_b(\r^2)$ and $x,y\in\r,$
	$$\ll|\bar P_tg(x,y) - P^N_tg(x,y)\rr|\leq C_t\cdot N^{-1/2}\cdot ||g||_{3,\infty} + C_t\cdot N^{-1}\ll(1+x^2\rr)\cdot ||g||_{2,\infty}.$$
\end{theorem}

\begin{figure}[!h]
\begin{center}
	\includegraphics[trim = 5.8cm 17cm 0cm 4cm,clip,width=0.95\linewidth]{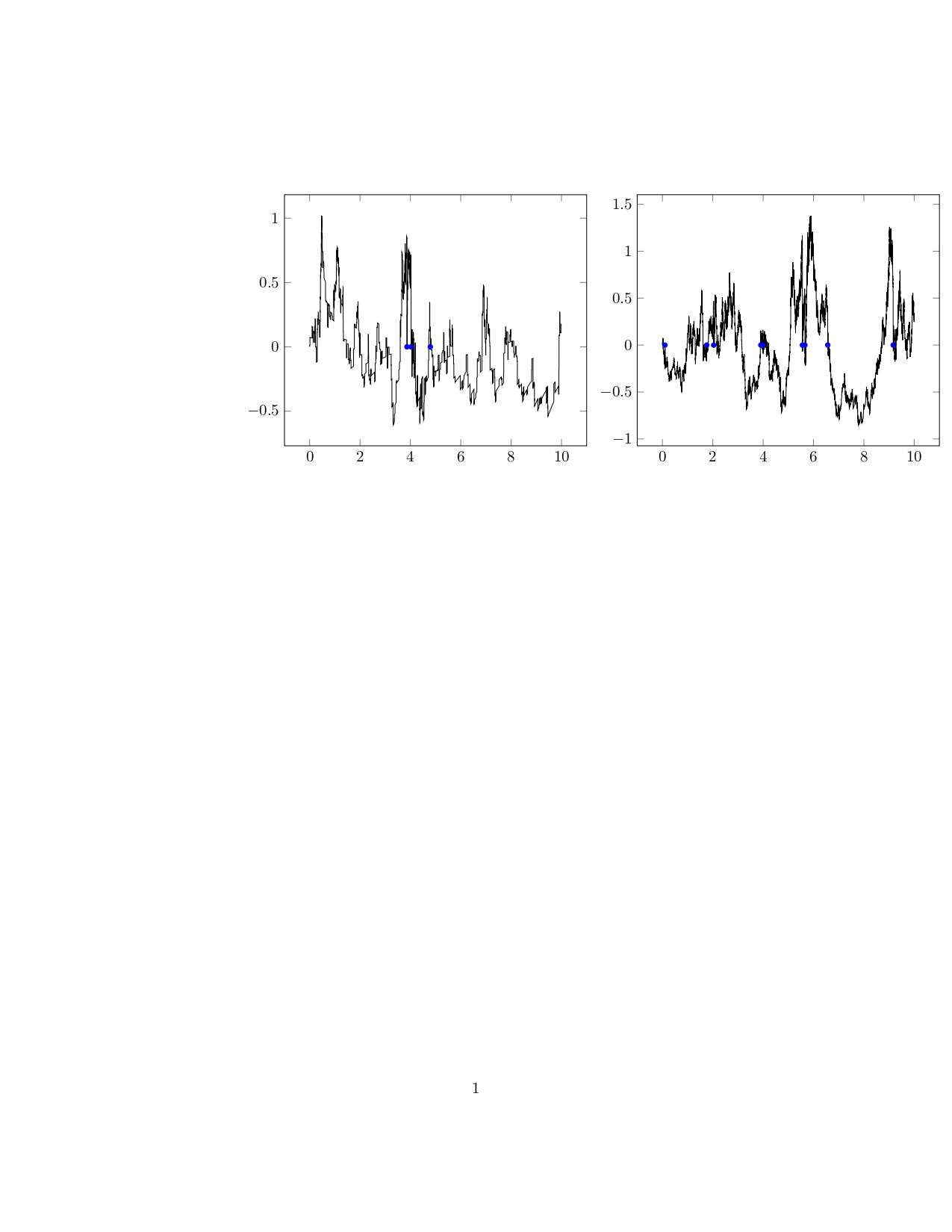}
	\caption{Simulation of trajectories of $(Y_t^{N})_{0\leq t\leq 10}$ with $X^{N}_0=0$, $Y^N_0=0$, $\alpha_1=\alpha_2=1$, $\nu=\mathcal{N}(0,1)$, $f_1(x)=f_2(x)=4/(1+\exp(2(1-x))),$ $N=100$ (left panel) and $N=1000$ (right panel). The blue points are the reset jump times of the process $Y^N.$}
	\label{simulationYN}
\end{center}
\end{figure}

Figure~\ref{simulationYN} shows two simulations of the process $Y^{N}$ with $N=100$ and $N=1000.$

\section{Proof of Theorem~\ref{mainresult}}\label{proof}

{We prove Theorem~\ref{mainresult} in three steps. In Section~\ref{generator}, we obtain an explicit convergence speed for the convergence of the generators $A^N$ of $(X^N,Y^N)$ to the generator $\bar A$ of $(\bar X,\bar Y)$ (see Lemma~\ref{convergencegenerator}). Then, we want to use this result using the following formula
	\begin{equation}\label{trotterkato}
	\bar P_tg(x,y) - P^N_tg(x,y) = \int_0^t P^N_{t-s}\ll(\bar A - A^N\rr)\bar P_sg(x,y)ds,
	\end{equation}
	to deduce the convergence speed for the semigroups $P^N$ and $\bar P.$ The problem is that our bound for the convergence speed of $A^Ng(x,y) - \bar Ag(x,y)$ depends on the derivatives of the function~$g$. So we need 	a control on the derivatives of the limit semigroup $(\bar P_t)_t$. This is done at Section~\ref{semigroups} (Proposition~\ref{semigroupc3}). Finally, at Section~\ref{proofmainresult}, we apply the formula~\eqref{trotterkato} to conclude the proof.}

\subsection{Convergence of the infinitesimal generators}
\label{generator}

One can note that the processes $(\bar X,\bar Y)$ and $(X^N,Y^N)$ ($N\in\n^*$) are Markov processes. In the following, we note $\bar P$ and $P^N$ ($N\in\n^*$) their respective semigroups, and $\bar A$ and $A^N$ ($N\in\n^*$) their respective infinitesimal generators. {Usually, if $Z$ is a $\r^2-$valued Markov process, we define its semigroup $P$ by: for any $t\geq 0,z\in\r^2,$ and any continuous and bounded function~$g$,
	$$P_tg(z) = \espcc{z}{g(Z_t)},$$
	where the index $z$ indicates that we consider the law under which $Z_0=z.$
	
	The definition of the semigroup of a Markov process coincides in general to the one above, but the notion of infinitesimal generators can change from a reference to another. Here we consider a notion of extended generators (similar as the one used in \cite{meyn_stability_1993, davis_markov_1993}). If $Z$ is a $\r^2-$valued Markov process (whose generator will be denoted by $A$), we note $\mathcal{D}'(A)$ the set of continuous and bounded functions~$g$ for which it exists a function $h_g$ that satisfies, for all $z\in\r^2, t\geq 0,$
	$$\mathbb{E}_z \int_0^t |h_g(Z_s)|ds<\infty~~\textrm{ and }~~\espcc{z}{g(Z_t)} - g(z) = \mathbb{E}_z \int_0^t h_g(Z_s)ds.$$
	If this function $h_g$ exists, we note it $Ag := h_g.$ Then $A$ is an operator called the (extended) generator of~$Z$. Note that the expression of a generator does not depend on the different definitions, what can change from a definition to another is the domain on which the generator is defined.
	
	\begin{example}
		Let $X$ be a Markov process solution to the following one-dimensional SDE
		$$dX_t = b(X_t)dt + \sigma(X_t)dW_t + \int_{\r_+}\Phi(X_{t-},u)\uno{z\leq f(X_{t-})}d\pi(t,z,u),$$
		where $W$ is a standard Brownian motion and $\pi$ a Poisson measure on $\r_+^2\times\r$ with intensity $dt\cdot dz\cdot d\nu(u).$ Then, Ito's formula allows to prove (provided some sufficient control on the functions $b,\sigma,f,\Phi$) that $\mathcal{D}'(A)$ contains $C^2_b(\r)$ and that
		$$Ag(x) = b(x)g'(x) + \frac12 \sigma(x)^2g''(x) + f(x)\int_\r \ll[g(x + \Phi(x,u))-g(x)\rr]d\nu(u).$$
		
		In this example we consider a one-dimensional SDE to keep things simple, but this can be generalized to any finite-dimensional SDE.
	\end{example}
}

We refer to Appendix~$A$ of \cite{erny_mean_2022} for the study of some basic properties of this notion of generators (in particular for Proposition~5.6 where~\eqref{trotterkato} is proved for another model).

Let us give the form of the generators $A^N$ and $\bar A$ of the processes $(X^N,Y^N)$ and $(\bar X,\bar Y).$ The following lemma is a direct consequence of Ito's formula and of the form of the SDEs~\eqref{XNYN} and~\eqref{barXbarY}.
\begin{lemma}\label{domaingenerator}
	Grant Assumptions~\ref{hyp}$.(i)-(ii)$.
	$\mathcal{D}'(A^N)$ contains $C^1_b(\r^2)$, and for $g\in C^1_b(\r^2)$ and $(x,y)\in\r^2,$
	\begin{align*}A^Ng(x,y) =& -\alpha_1 x\partial_1g(x,y) - \alpha_2 y\partial_2g(x,y) \\
	&+ Nf_1(x)\int_\r\ll[g\ll(x+\frac{u}{\sqrt{N}} - \frac{x}{N},y+\frac{u}{\sqrt{N}}\rr)-g(x,y)\rr]d\nu(u) \\
	&+ f_2(y)\ll[g(x,0)-g(x,y)\rr].
	\end{align*}
	
	$\mathcal{D}'(\bar A)$ contains $C^2_b(\r^2)$, and for $g\in C^2_b(\r^2)$ and $(x,y)\in\r^2,$
	\begin{align*}
	\bar Ag(x,y) =& -\alpha_1 x\partial_1g(x,y) - x f_1(x)\partial_1g(x,y) - \alpha_2 y\partial_2g(x,y)\\
	& + \frac{\sigma^2}{2}f_1(x)\sum_{i,j=1}^2\partial^2_{i,j}g(x,y) + f_2(y)\ll[g(x,0)-g(x,y)\rr].
	\end{align*}
\end{lemma}

We have the following control on the convergence of the generators.
\begin{lemma}\label{convergencegenerator}
	Grant Assumptions~\ref{hyp}$.(i)-(ii)$ and Assumption~\ref{hyp}$.(iv)$. For $g\in C^3_b(\r^2)$ and $(x,y)\in\r^2,$
	$$\ll|A^Ng(x,y) - \bar Ag(x,y)\rr|\leq \frac{f_1(x)}{6 N^{1/2}}\ll(\int_\r |u|^3d\nu(u)\rr) \sum_{|\beta|=3}||\partial_\beta g||_\infty + \frac12 x^2f_1(x) N^{-1}\ll|\ll|\partial^2_{11}g\rr|\rr|_\infty.$$
\end{lemma}

\begin{proof}
	
	Recalling that $\nu$ is a centered probability measure and that $\sigma^2$ is its variance,
	\begin{align*}
	&A^Ng(x,y) - \bar Ag(x,y)\\
	=& \int_\r Nf_1(x)\ll[g\ll(x+\frac{u}{\sqrt{N}}- \frac{x}{N},y+\frac{u}{\sqrt{N}}\rr)-g(x,y)\rr]d\nu(u) \\
	&+ xf_1(x)\partial_1g(x,y) - \frac{\sigma^2}{2}f_1(x)\sum_{i,j=1}^2\partial^2_{i,j}g(x,y)\\
	=&Nf_1(x)\int_\r\ll[g\ll(x+\frac{u}{\sqrt{N}}-\frac{x}{N},y+\frac{u}{\sqrt{N}}\rr)-g(x,y) - \frac{u}{\sqrt{N}}\sum_{i=1}^2\partial_ig(x,y) + \frac{x}{N}\partial_1g(x,y)\right.\\
	&\hspace*{2.2cm}\ll.- \frac{u^2}{2N}\sum_{i,j=1}^2\partial^2_{i,j}g(x,y)\rr]d\nu(u).
	\end{align*}
	
	In order to make appear the second-order term of the Taylor polynomial of $g(x+u/\sqrt{N} - x/N,y+u/\sqrt{N}),$ let us write
	\begin{multline*}
	\frac{u^2}{2N}\sum_{i,j=1}^2\partial^2_{i,j}g(x,y)\\
	=\frac12\ll(\frac{u}{\sqrt{N}} - \frac{x}{N}\rr)^2\partial^2_{11}g(x,y) + \frac{u^2}{2N}\partial^2_{22}g(x,y) + \frac{u}{\sqrt{N}}\ll(\frac{u}{\sqrt{N}} - \frac{x}{N}\rr)\partial^2_{12}g(x,y)\\
	+\frac{ux}{N\sqrt{N}}\partial^2_{11}g(x,y) - \frac{x^2}{2N^2}\partial^2_{11}g(x,y) + \frac{ux}{N\sqrt{N}}\partial^2_{12}g(x,y).
	\end{multline*}
	Therefore, since $\int_\r ud\nu(u)=0,$
	\begin{align*}
	&A^Ng(x,y) - \bar Ag(x,y)\\
	=&Nf_1(x)\int_\r\ll[g\ll(x+\frac{u}{\sqrt{N}}-\frac{x}{N},y+\frac{u}{\sqrt{N}}\rr)-g(x,y) - \frac{u}{\sqrt{N}}\sum_{i=1}^2\partial_ig(x,y) + \frac{x}{N}\partial_1g(x,y)\right.\\
	&\ll.- \ll(\frac12\ll(\frac{u}{\sqrt{N}} - \frac{x}{N}\rr)^2\partial^2_{11}g(x,y) + \frac{u^2}{2N}\partial^2_{22}g(x,y) + \frac{u}{\sqrt{N}}\ll(\frac{u}{\sqrt{N}} - \frac{x}{N}\rr)\partial^2_{12}g(x,y)\rr)\rr]d\nu(u)\\
	&+ \frac{x^2}{2N}f_1(x)\partial^2_{11}g(x,y).
	\end{align*}
	
	As a consequence
	\begin{multline*}
	\ll|A^Ng(x,y) - \bar Ag(x,y)\rr|\\
	\leq Nf_1(x)\int_\r\ll[\ll|g\ll(x+\frac{u}{\sqrt{N}} - \frac{x}{N},y+\frac{u}{\sqrt{N}}\rr)-g(x,y) \rr.\rr.\\
	- \frac{u}{\sqrt{N}}\sum_{i=1}^2\partial_ig(x,y) + \frac{x}{N}f_1(x)\partial_1g(x,y)
	 - \frac12\ll(\frac{u}{\sqrt{N}} 
	 - \frac{x}{N}\rr)^2\partial^2_{11}g(x,y) \\
	\ll.\ll.- \frac{u^2}{2N}\partial^2_{22}g(x,y) - \frac{u}{\sqrt{N}}\ll(\frac{u}{\sqrt{N}} - \frac{x}{N}\rr)\partial^2_{12}g(x,y)\rr|\rr]d\nu(u)\\
	+  \frac{x^2}{2N}f_1(x)|\partial^2_{11}g(x,y)|.
	\end{multline*}
	Then Taylor-Lagrange's inequality gives the result.
\end{proof}

\subsection{Regularity of the limit semigroup and stochastic flow}
\label{semigroups}

In this section, we prove the regularity of the semigroup of the limit process $(\bar X,\bar Y)$ w.r.t. its initial condition and a control of its derivatives. However, there are additional technical difficulties in this proof compared to \cite{erny_mean_2022}. {We rely on Theorem~1.4.1 of \cite{kunita_lectures_1986} which guarantees the differentiability of stochastic flows under some conditions. Let us recall that the stochastic flow of $(\bar X, \bar Y)$ is the family of functions $(x,y)\mapsto (\bar X^{(x)}_t,\bar Y^{(x,y)}_t)$ ($t\geq 0$). A well-known kind of flow is called "Brownian flow". By definition, a Brownian flow is a stochastic flow that is almost surely continuous w.r.t. the time-parameter~$t$ with independent increments w.r.t. the time-parameter.}

The main difficulty is that the trajectories of the process $(\bar X, \bar Y)$ are not continuous. Hence, the stochastic flow of this process is not a Brownian flow and Theorem~1.4.1 of \cite{kunita_lectures_1986} cannot be applied. So, we study the stochastic flow between its jump times. Then we face another difficulty: the jump times of the stochastic flow depends on the initial condition, and we do not want to study the regularity of this dependency.

This is why we introduce the following auxiliary limit process $(\tilde X,\tilde Y)$. The only difference with the limit process $(\bar X,\bar Y)$ is that the stochastic intensity of the jump term in the SDE of $\tilde Y$ is constant. In other words, the occurrences of the jump times of the process $\tilde Y$ are distributed as a homogeneous Poisson process.
\begin{align*}
\tilde X_t=&\tilde X_0 -\alpha_1\int_0^t\tilde X_sds + \sigma\int_0^t\sqrt{f_1(\tilde X_s)}dW_s - \int_0^t\tilde X_s f_1(\tilde X_s)ds,\\
\tilde Y_t=&\tilde Y_0 - \alpha_2\int_0^t\tilde Y_sds + \sigma\int_0^t\sqrt{f_1(\tilde X_s)}dW_s - \int_{[0,t]\times [0,||f_2||_\infty]} \tilde Y_{s-}d\tilde\pi^1(s,z),
\end{align*}
where $\tilde \pi^1$ is a Poisson measure on~$\r_+^2$ with Lebesgue intensity.

Let us note $N_t = \tilde\pi^1([0,t]\times||f_2||_\infty)$ the number of jumps of the process~$\tilde Y$ before time~$t.$ By definition of $\tilde \pi^1,$ $(N_t)_{t\geq 0}$ is a homogeneous Poisson process with rate $||f_2||_\infty.$ We also note $T_n$ ($n\geq 1$) the jump times of the process~$\tilde Y.$

The next proposition is a consequence of Girsanov's theorem for jump processes. It gives the Radon-Nikodym derivative between the distribution of~$(\bar X_t,\bar Y_t)$ and that of~$(\tilde X_t,\tilde Y_t),$ for any $t\geq 0.$ This result is a straightforward application of Theorem~$VI.T3$ of \cite{bremaud_point_1981} (see also of Theorem~3.5 \cite{locherbach_likelihood_2002}), its proof is therefore omitted.
\begin{proposition}\label{girsanov}
	Grant Assumptions~\ref{hyp}$.(i)-(ii)$. Let $N_t:=\tilde\pi^1([0,t]\times[0,||f_2||_\infty])$ and $T_n$ ($n\in\n^*$) be the jump times of $(N_t)_t.$ For all $g\in C_b(\r^2),$ and, $x,y\in\r,$
	\begin{multline*}
	\esp{g\ll(\bar{X}^{(x)}_t,\bar{Y}^{(x,y)}_t\rr)}\\=\esp{g\ll(\tilde X^{(x)}_t,\tilde Y^{(x,y)}_t\rr)\prod_{n=1}^{N_t}\frac{f_2\ll(\tilde Y^{(x,y)}_{T_n}\rr)}{||f_2||_\infty}\exp\ll[{-\int_0^t\ll(f_2\ll(\tilde Y^{(x,y)}_s\rr)-||f_2||_\infty\rr)ds}\rr]}.
	\end{multline*}
\end{proposition}

In order to obtain a control on the regularity of the semigroup $\bar P$ of the limit process $(\bar X,\bar Y),$ we rely on the previous proposition and on a control of the derivatives of the stochastic flow of the auxiliary process $(\tilde X,\tilde Y).$ The next step is to prove this control.

Let us note $\phi_t(x,y)$ ($t\in\r_+,(x,y)\in\r^2$) the stochastic flow related to the following $2-$dimensional SDE:
\begin{align}
X_t=&X_0 -\alpha_1 \int_0^t X_sds + \sigma\int_0^t\sqrt{f_1(X_s)}dW_s - \int_0^t X_sf_1(X_s)ds,\label{XY}\\
Y_t=& Y_0 - \alpha_2 \int_0^tY_sds + \sigma\int_0^t\sqrt{f_1(X_s)}dW_s.\nonumber
\end{align}
{With our previous notation, it means that $\phi_t(x,y) := (X^{(x)}_t,Y^{(x,y)}_t),$ if $(X,Y)$ is defined as the solution of the above $2-$dimensional SDE.} One can note that $\phi_t(x,y)$ is the flow of the process $(\tilde X,\tilde Y)$ between its jump times. {It is clear that $\phi_t(x,y)$ is a stochastic flow by Definition~1.1.1 of \cite{kunita_lectures_1986}. It is even a Brownian flow since the process $(X,Y)$ has independent increments and continuous trajectories.}

\begin{remark}
	Contrarily to the more general framework of \cite{kunita_lectures_1986}, the stochastic flow $\phi_t(x,y)$ depends only  on one time parameter, because the process is time-homogeneous.
\end{remark}

\begin{proposition}\label{flotc3}
	Grant Assumption~\ref{hyp}$.(i)$ and Assumption~\ref{hyp}$.(iii)$. Almost surely, for all $t\geq 0,$ $(x,y)\mapsto\phi_t(x,y)$ is a $C^3-$diffeomorphism.
\end{proposition}

\begin{proof}
	We rely on Theorem~1.4.1 of \cite{kunita_lectures_1986} and use the notation therein. {Let us note $(a,b)$ the local characteristics of the flow $\phi_t(x,y).$ The notion of local characteristics is formally defined in Assumption~1 of \cite[p. 8]{kunita_lectures_1986}, and when we study the flow related to an SDE with a drift term and a Brownian term, the function~$b$ can be shown to be the coefficient of the drift term, and $a$ the matrix of quadratic covariations of the Brownian term. Here the local characteristics of $\phi_t(x,y)$ are given by}
	$$b(x,y) = \begin{pmatrix} -\alpha_1 x - xf_1(x)\\-\alpha_2 y\end{pmatrix}\textrm{ and }a((x,y),(x',y'))=\sigma^2\sqrt{f_1(x)f_1(x')}\begin{pmatrix}1&1\\1&1\end{pmatrix}.$$
	
	One can check that the Assumptions~1,~2 and~3 of \cite[p. 8, 9, 15]{kunita_lectures_1986} are satisfied. {Indeed, Assumptions~1 and~2 follows easily from the fact that, if $(X,Y)$ is solution to~\eqref{XY}, then 
		$$\espcc{(x,y)}{(X_t)^2 + (Y_t)^2}\leq C_t(1+x^2+y^2),$$
		where the inequality above can be proved with the same reasoning as in the proof of Lemma~\ref{aprioriestimate}. For Assumption~3, it is sufficient to remark that the functions $x\mapsto xf_1(x)$ and $x\mapsto \sqrt{f_1(x)}$ are Lipschitz continuous under Assumptions~\ref{hyp}$.(i)$ and~\ref{hyp}$.(iii)$.}
	
	Then, Theorem~1.4.1 of \cite{kunita_lectures_1986} allows to conclude the proof.
\end{proof}

\begin{proposition}\label{flotregulier}
	Grant Assumptions~\ref{hyp}$.(i)-(iii)$. Almost surely, for all $t\geq 0,$ $(x,y)\mapsto (\tilde X_t^{(x)},\tilde Y_t^{(x,y)})$ is $C^3$ and, for all even $p\in\n^*,$ for all multi-index $\beta\in\n^2$ such that $1\leq |\beta|\leq 3,$
	$$\underset{x\in\r}{\sup}~\esp{\underset{0\leq s\leq t}{\sup}\ll(\partial_\beta \tilde X_s^{(x)}\rr)^p} + \underset{x,y\in\r}{\sup}\esp{\underset{0\leq s\leq t}{\sup}\left(\partial_\beta\tilde Y_s^{(x,y)}\right)^p}<\infty.$$
\end{proposition}

\begin{proof}
	We already know that, for all $t\geq 0,$ $x\mapsto\tilde X_t^{(x)}$ is $C^3$ by Proposition~\ref{flotc3}. Indeed, by definition $\tilde X_t^{(x)}$ is the first coordinate of $\phi_t(x,y).$
	
	Now, let us recall that $T_n$ ($n\in\n^*$) denote the jump times of the process $(\tilde X^{(x)},\tilde Y^{(x,y)})$ and that these jump times do not depend on the initial condition~$(x,y)$. For $t<T_1,$
	$$(\tilde X^{(x)}_t,\tilde Y^{(x,y)}_t) = \phi_t(x,y)$$
	is $C^3$ w.r.t. $(x,y)$ by Proposition~\ref{flotc3}. And, for $t\in[T_n,T_{n+1}[$ ($n\in\n^*$),
	$$(\tilde X^{(x)}_t,\tilde Y^{(x,y)}_t) = \phi_{t-T_n}(\tilde X^{(x)}_{T_n},0),$$
	which is also $C^3$ w.r.t. $(x,y)$ as a composition of the function $(x,y)\in\r^2\mapsto\phi_{t-T_n}(x,y)$ and the function~$(x,y)\in\r^2\mapsto(\tilde X^{(x)}_{T_n},0).$ This proves that, almost surely, for all $t\geq 0,$ $(x,y)\mapsto(\tilde X_t^{(x)},\tilde Y^{(x,y)}_t)$ is~$C^3.$ Now, we prove the second part of the statement.
	
	By definition, for all $t\geq 0,$
	$$\tilde X^{(x)}_t = x  - \alpha_1\int_0^t \tilde X^{(x)}_sds - \int_0^t \tilde X_s^{(x)}f_1(\tilde X_s^{(x)})ds + \sigma\int_0^t\sqrt{f_1\ll(\tilde X^{(x)}_s\rr)}dW_s.$$
	
	This implies
	\begin{align*}
	\partial_x\tilde X^{(x)}_t =& 1  - \alpha_1\int_0^t \partial_x\tilde X^{(x)}_sds - \int_0^t \tilde X_s^{(x)}\partial_x \tilde X_s^{(x)}f_1'(\tilde X_s^{(x)})ds- \int_0^t \partial_x\tilde X_s^{(x)}f_1(\tilde X_s^{(x)})ds\\
	&+ \sigma\int_0^t\partial_x\tilde X^{(x)}_s\ll(\sqrt{f_1}\rr)'\ll(\tilde X^{(x)}_s\rr)dW_s.
	\end{align*}
	
	We also know that
	\begin{multline*}\partial_x\tilde Y^{(x,y)}_t = -\alpha_2\int_0^t\partial_x\tilde Y^{(x,y)}_sds + \sigma\int_0^t\partial_x\tilde X^{(x)}_s\cdot \sqrt{f_1\ll(\tilde X^{(x)}_s\rr)}dW_s \\
	- \int_0^t\int_{[0,||f_2||_\infty]}\partial_x\tilde Y^{(x,y)}_sd\tilde\pi(s,z).
	\end{multline*}
	
	Noticing that the initial condition of the process $(\partial_x\tilde X^{(x)}_t,\partial_x \tilde Y^{(x,y)}_t)$ is~$(1,0)$, Lemma~\ref{aprioriestimate} allows to prove the result of the lemma for the first order partial derivative w.r.t.~$x.$ The results for the other partial derivatives follow from the same reasoning.
\end{proof}

\begin{proposition}\label{semigroupc3}
	Grant Assumptions~\ref{hyp}$.(i)-(iii)$. For all $t\geq 0$ and $g\in C^3_b(\r^2),$ the function $(x,y)\in\r^2\mapsto \bar P_tg(x,y)$ is $C^3$ and for all multi-index $\beta$ such that $|\beta|=3,$
	$$\underset{0\leq s\leq t}{\sup}\ll|\ll|\partial_\beta\bar P_sg\rr|\rr|_\infty \leq C_t||g||_{3,\infty},$$
	where $C_t>0$ does not depend on~$g.$
\end{proposition}

\begin{proof}
	For the sake of readability, we do not prove the result for the multi-indexes $\beta$ such that $|\beta|=3$ but for those such that $|\beta|=1.$ This case is easier to prove and the techniques are exactly the same.
	
	Recalling Proposition~\ref{girsanov}, we have
	\begin{equation}
	\label{barPtg}
	\bar P_tg(x,y)=\esp{g\ll(\tilde X^{(x)}_t,\tilde Y^{(x,y)}_t\rr)\prod_{n=1}^{N_t}\frac{f_2\ll(\tilde Y^{(x,y)}_{T_n}\rr)}{||f_2||_\infty}\exp\ll[{-\int_0^t\ll(f_2\ll(\tilde Y^{(x,y)}_s\rr)-||f_2||_\infty\rr)ds}\rr]}.
	\end{equation}
	
	To prove that the function $\bar P_tg$ is $C^1,$ we rely on Lemma~\ref{derivevitali} which is a generalization of the classical result about the "differentiation under the integral sign" that uses Vitali theorem instead of the dominated convergence theorem.
	
	By Proposition~\ref{flotregulier}, we know that the expression in the expectation of the right-hand side of~\eqref{barPtg} is $C^1$ almost surely. In addition, one can control the expectation of the square of the derivatives of this expression thanks to Proposition~\ref{flotregulier} (see below for the explicit expression of the derivatives). More precisely, if we note $Z^{(x,y)}$ the random variable appearing in the expectation of~\eqref{barPtg}, then Proposition~\ref{flotregulier} implies that for any $x',y'\in\r,$
		$$\underset{x\in\r}{\sup} \int_\Omega \ll|\partial_x Z^{(x,y')}\rr|^2d\mathbb{P}(\omega)<\infty\textrm{ and }\underset{y\in\r}{\sup} \int_\Omega \ll|\partial_y Z^{(x',y)}\rr|^2d\mathbb{P}(\omega)<\infty.$$
		This implies that the uniform integrability condition of Lemma~\ref{derivevitali} is satisfied.
	
	So Lemma~\ref{derivevitali} (applied twice since there are two coordinates) implies that $\bar P_tg$ is~$C^1$. In addition, denoting by
	$$L_t(x,y) := \prod_{n=1}^{N_t}\frac{f_2\ll(\tilde Y^{(x,y)}_{T_n}\rr)}{||f_2||_\infty}\exp\ll[{-\int_0^t\ll(f_2\ll(\tilde Y^{(x,y)}_s\rr)-||f_2||_\infty\rr)ds}\rr],$$
	we have {(thanks to Lemma~\ref{derivevitali})}
	\begin{align*}
	&\partial_2\ll(\bar P_tg\rr)(x,y)=\esp{\partial_y \tilde Y^{(x,y)}_t\cdot\partial_2g\ll(\tilde X^{(x)}_t,\tilde Y^{(x,y)}_t\rr)L_t(x,y)}\\
	&+\esp{g\ll(\tilde X^{(x)}_t,\tilde Y^{(x,y)}_t\rr)\sum_{n=1}^{N_t}\partial_y\tilde Y^{(x,y)}_{T_n}\frac{f_2'\ll(\tilde Y^{(x,y)}_{T_n}\rr)}{||f_2||_{\infty}}\prod_{k\neq n}\frac{f_2\ll(\tilde Y^{(x,y)}_{T_k}\rr)}{||f_2||_\infty}\rr.\\
&\ll.\quad\quad\quad\times		
	\exp\ll(\int_0^t(||f_2||_\infty-f_2(\tilde Y^{(x,y)}_s))ds\rr)}\\
	&-\esp{\int_0^t\partial_y\tilde Y^{(x,y)}_sf_2'\ll(\tilde Y^{(x,y)}_s\rr)ds \cdot g\ll(\tilde X^{(x)}_t,\tilde Y^{(x,y)}_t\rr)L_t(x,y)}.
	\end{align*}
	
	As a consequence,
	\begin{align*}
	\ll|\partial_2\ll(\bar P_tg\rr)(x,y)\rr|\leq& ||\partial_2g||_\infty e^{t||f_2||_\infty} \esp{\ll|\partial_y\tilde Y^{(x,y)}_t\rr|}\\
	&+ ||g||_\infty \frac{||f'_2||_\infty}{||f_2||_\infty} e^{||f_2||_\infty t} \esp{N_t\underset{0\leq s\leq t}{\sup}\ll|\partial_y\tilde Y^{(x,y)}_{s}\rr|}\\
	&+ ||f_2'||_\infty\cdot||g||_\infty e^{t||f_2||_\infty}t\esp{\underset{0\leq s\leq t}{\sup}\ll|\partial_y\tilde Y^{(x,y)}_s\rr|}.
	\end{align*}
	
	The expressions of the first line and the third line above can be bounded uniformly in~$(x,y)$ thanks to Proposition~\ref{flotregulier}. For the expression of the second line, we can use Cauchy-Schwarz's inequality for example (we can recall that $N_t$ follows the Poisson distribution with parameter $t||f_2||_\infty$).
	
	This implies that
	$$\underset{x,y\in\r}{\sup}~\underset{0\leq s\leq t}{\sup} \ll|\partial_2\ll(\bar P_sg\rr)(x,y)\rr|\leq C_t||g||_{1,\infty}.$$
	
	With similar computation, we can prove that
	$$\underset{x,y\in\r}{\sup}~\underset{0\leq s\leq t}{\sup} \ll|\partial_1\ll(\bar P_sg\rr)(x,y)\rr|\leq C_t||g||_{1,\infty}.$$
	
	So we have proved the result for the two multi-indexes $\beta$ such that $|\beta|=1.$ Note that in the case for the multi-indexes $\beta$ such that $|\beta|=3,$ we would need a control on higher moments of the processes $\tilde X,$ $\tilde Y$ and their derivatives to control expectations of product of these processes. Indeed, the partial derivative w.r.t. the multi-index~$(1,2)$ would involve terms as
	$$\esp{\ll(\underset{0\leq s\leq t}{\sup}|\partial_x \tilde Y^{(x,y)}_s|\rr)\ll(\underset{0\leq s\leq t}{\sup}|\partial_y \tilde Y^{(x,y)}_s|\rr)^2}$$
	which can be controlled using Cauchy-Schwarz's inequality, recalling that Proposition~\ref{flotregulier} allows to bound every polynomial moment of the partial derivatives of the process $(\tilde X,\tilde Y)$ w.r.t. its initial condition.
\end{proof}

\subsection{End of the proof of Theorem~\ref{mainresult}}
\label{proofmainresult}

With a similar proof to the one of Proposition~5.6 of \cite{erny_mean_2022}, we can prove that, for all $t\geq 0,$ $N\in\n^*,$ $g\in C^3_b(\r^2),$ $x,y\in\r,$
$$\ll(\bar P_t - P^N_t\rr)g(x,y) = \int_0^t P^N_{t-s}\ll(\bar A-A^N\rr)\bar P_{s}g(x,y)ds.$$

Consequently,
\begin{align*}
\ll|\ll(\bar P_t - P^N_t\rr)g(x,y)\rr|\leq& \int_0^t\espcc{(x,y)}{\ll|\ll(\bar A-A^N\rr)\bar P_sg(X^N_{t-s},Y^N_{t-s})\rr|}ds\\
\leq& C\cdot N^{-1/2}\int_{0}^t\espcc{(x,y)}{f_1(X^N_{t-s})\sum_{|\beta|=3}||\partial_\beta\ll(\bar P_sg\rr)||_\infty}ds \\
&+ C N^{-1}\int_0^t \espcc{(x,y)}{\ll(X^N_{t-s}\rr)^2f_1\ll(X^N_{t-s}\rr)||\partial^2_{11}(\bar P_sg)||_\infty}ds\\
\leq& C_t\cdot N^{-1/2}\cdot ||f_1||_\infty ||g||_{3,\infty} + C_t\cdot N^{-1}\cdot ||f_1||_\infty\ll(1+x^2\rr)||g||_{2,\infty},	
\end{align*}
where we have used Lemma~\ref{convergencegenerator} to obtain the second inequality above, and Proposition~\ref{semigroupc3} and Lemma~\ref{aprioriXNYN} to obtain the last one. {More precisely, Proposition~\ref{semigroupc3} is used to control the partial derivatives of the function $\bar P_s g,$ and Lemma~\ref{aprioriXNYN} to control the second order moment of $X^N_{t-s}.$}

\section{Technical lemmas}

The first lemma allows to differentiate under an integral sign. Contrary to the classical result that relies on the dominated convergence theorem, this result relies on Vitali convergence theorem. This result is a direct consequence of Lemma~6.1 of \cite{eldredge_strong_2018}.

\begin{lemma}\label{derivevitali}
	Let $F:(x,\omega)\in\r\times\Omega\mapsto\r$ be a measurable function, where $(\Omega,\mathbb{P})$ is a probability space. Assume that:
	\begin{itemize}
		\item $\mathbb{P}-$almost surely, $x\mapsto F(x,\omega)$ is~$C^1$,
		\item the set of functions $\{\omega\mapsto \partial_1 F(x,\omega)~:~x\in\r\}$ is uniformly integrable.
	\end{itemize}
	
	Then, $\mathbb{P}-$almost surely, $x\mapsto \int_\Omega F(x,\omega)d\mathbb{P}(\omega)$ is $C^1$ and its derivative is
	$$\int_\Omega\partial_1 F(x,\omega)d\mathbb{P}(\omega).$$
\end{lemma}

\begin{proof}

	By Lemma~6.1 of \cite{eldredge_strong_2018}, for any $T>0,$ the function $x\mapsto \int_\Omega F(x,\omega)d\mathbb{P}(\omega)$ is almost surely $C^1$ on $[-T,T]$ with derivative
		$$\int_\Omega \partial_1 F(x,\omega)d\mathbb{P}(\omega).$$
		As it holds true for any $T>0,$ the result is proved.
\end{proof}

The next lemma allows to obtain a priori estimates on solutions of certain SDEs (in our model, it is used to control the polynomial moments of the derivatives of the limit process in Proposition~\ref{flotregulier}). This result is not written to be general, because some specific properties of our SDEs allow some simplifications to obtain the result.

\begin{lemma}\label{aprioriestimate}
	Let $(X_t,Y_t)_t$ be a solution of the following two-dimensional SDE
	\begin{align*}
	dX_t =& b_1(X_t,Y_t)dt + \varsigma_1(X_t,Y_t)dW_t,\\
	dY_t =& b_2(X_t,Y_t)dt + \varsigma_2(X_t,Y_t)dW_t - \int_{[0,+\infty[}Y_{t-}\uno{z\leq f(Y_{t-})}d\pi(z,t),
	\end{align*}
	where $(W_t)_t$ is one-dimensional standard Brownian motion, $\pi$ a Poisson measure on $\r_+^2$ with Lebesgue intensity and $f$ a non-negative measurable function. Assume that the functions $b_1,$ $b_2,$ $\varsigma_1$ and $\varsigma_2$ are sublinear: for all $x,y\in\r,$
	$$|b_1(x,y)| + |b_2(x,y)| + |\varsigma_1(x,y)| + |\varsigma_2(x,y)|\leq C(1+|x|+|y|).$$
	
	Then:
	\begin{itemize}
		\item[$(i)$] for any even $p\in\n$ such that $X_0$ and $Y_0$ both have a finite $p-$order moment, for all $t\geq 0,$
		$$\esp{X_t^p}+\esp{Y_t^p}\leq C_{t,p}(1+\esp{X_0^p} + \esp{Y_0^p}),$$
		with $C_{t,p}>0$ some constant,
		\item[$(ii)$] if $X_0$ and $Y_0$ both have a finite $2p-$th order moment (for some even $p\in\n^*$), for all $t\geq 0,$
		$$\esp{\underset{0\leq s\leq t}{\sup}X_s^p} + \esp{\underset{0\leq s\leq t}{\sup}Y_s^p}\leq C_{t,p}\ll(1+\esp{X_0^{2p}} + \esp{Y_0^{2p}}\rr),$$
		with $C_{t,p}>0$ another constant.
	\end{itemize}
\end{lemma}

\begin{proof}
	{\it Step~1.} Let us prove the point~$(i)$. By Ito's formula, for all $t\geq 0,$
	\begin{align}
	X_t^p=& X_0^p  +p\int_0^t X_s^{p-1} b_1(X_s,Y_s)ds +\frac{p(p-1)}{2}\int_0^t X_s^{p-2}\varsigma_1(X_s,Y_s)^2ds \nonumber\\&+ p\int_0^t X_s^{p-1}\varsigma_1(X_s,Y_s)dW_s,\label{eqx}\\
	Y_t^p=& Y_0^p  +p\int_0^t Y_s^{p-1} b_2(X_s,Y_s)ds +\frac{p(p-1)}{2}\int_0^t Y_s^{p-2}\varsigma_2(X_s,Y_s)^2ds \nonumber\\
	& + p\int_0^t Y_s^{p-1}\varsigma_2(X_s,Y_s)dW_s - \int_{[0,t]\times[0,+\infty[} Y_{s-}^p\uno{z\leq f(Y_{s-})}d\pi(s,z)\nonumber\\
	\leq & Y_0^p  +p\int_0^t Y_s^{p-1} b_2(X_s,Y_s)ds +\frac{p(p-1)}{2}\int_0^t Y_s^{p-2}\varsigma_2(X_s,Y_s)^2ds \nonumber\\
	&+ p\int_0^t Y_s^{p-1}\varsigma_2(X_s,Y_s)dW_s.\label{eqy}
	\end{align}
	
	Now for $M>0,$ let us introduce the following stopping time
	$$\tau_M := \inf\{t>0~:~|X_t|>M\textrm{ or }|Y_t|>M\},$$
	
	and the function
	$$u_{M,p}(t) := \esp{(X_{t\wedge\tau_M})^p + (Y_{t\wedge\tau_M})^p}.$$
	
	Thanks to~\eqref{eqx}, we have for all $t\geq 0,$
	\begin{align*}
	&\esp{(X_{t\wedge\tau_M})^p}= \esp{X_0^p} + \esp{\int_0^{t\wedge\tau_M} \ll(pX_s^{p-1} b_1(X_s,Y_s) + \frac{p(p-1)}{2}X_s^{p-2}\varsigma_1(X_s,Y_s)^2\rr)ds}\\
	&\leq  \esp{X_0^p} + C\int_0^{t} \esp{|X_{s\wedge\tau_M}|^{p-1}\cdot |b_1(X_{s\wedge\tau_M},Y_{s\wedge\tau_M})| + X_{s\wedge\tau_M}^{p-2}\varsigma_1(X_{s\wedge\tau_M},Y_{s\wedge\tau_M})^2}ds\\
	&\leq \esp{X_0^p} + C\int_0^t\ll(1+u_{M,p}(s)\rr)ds,
	\end{align*}
	where we have used the sublinearity of the functions $b_1$ and $\varsigma_1$ to obtain the last line. With the same reasoning, we also have using~\eqref{eqy} for all $t\geq 0,$
	$$\esp{(Y_{t\wedge\tau_M})^p}\leq \esp{Y_0^p} + C\int_0^t\ll(1+u_{M,p}(s)\rr)ds.$$
	
	So, for all $t\geq 0,$
	$$u_{M,p}(t)\leq \esp{X_0^p} + \esp{Y_0^p} + Ct + C\int_0^tu_{M,p}(s)ds.$$
	
	Then, as the function $u_{M,p}$ is locally integrable (thanks to the stopping time $\tau_M$), Gr\"onwall's lemma implies that for all $t\geq 0,$
	\begin{equation}
	\label{umt}
	u_{M,p}(t)\leq \ll(\esp{X_0^p} + \esp{Y_0^p} + Ct\rr)e^{Ct}\leq C_t(1+\esp{X_0^p} + \esp{Y_0^p}),
	\end{equation}
	where $C_t>0$ does not depend on~$M.$
	
	Now let us prove that $\tau_M$ goes to infinity almost surely as $M$ goes to infinity. Firstly, as $(\tau_M)_M$ is non-decreasing w.r.t. $M$, we know that $\tau_M$ has an almost sure limit $\tau$ as $M$ goes to infinity. We have, for all $t>0,$ by Markov's inequality and thanks to~\eqref{umt},
	$$\pro{\tau_M\leq t}\leq \pro{X_{t\wedge\tau_M}^p+Y_{t\wedge\tau_M}^p\geq M^p}\leq \frac{1}{M^p}u_{M,p}(t)\underset{M\rightarrow\infty}{\longrightarrow}0.$$
	
	Consequently, for all $t>0,$
	$$\pro{\tau\leq t}\leq \pro{\bigcap_{M\in\n^*} \{\tau_M\leq t\}} = \underset{M\rightarrow\infty}{\lim} \pro{\tau_M\leq t} = 0.$$
	
	Hence,
	$$\pro{\tau<\infty}\leq \pro{\bigcup_{t\in\n^*}\{\tau\leq t\}}=0.$$
	This means that $\tau_M$ goes to infinity almost surely as $M$ goes to infinity.
	
	Finally, applying Fatou's lemma to~\eqref{umt} gives that for all $t\geq 0,$
	$$\esp{X_t^p+Y_t^p}\leq C_t(1+\esp{X_0^p} + \esp{Y_0^p}).$$
	
	{\it Step~2.} Let us prove the point~$(ii)$. By~\eqref{eqx} and Burkh\"older-Davis-Gundy's inequality,
	\begin{multline*}
	\esp{\underset{0\leq s\leq t}{\sup}X_s^p}
	\leq \esp{X_0}^p + p\int_0^t\esp{|X_s^{p-1}b_1(X_s,Y_s)|}ds \\
	+ \frac{p(p-1)}{2}\int_0^t\esp{X_s^{p-2}\varsigma_1(X_s,Y_s)^2}ds+C p\ll(\int_0^t\esp{X_s^{2p-2}\varsigma_1(X_s,Y_s)^2}ds\rr)^{1/2}\\
	\leq \esp{X_0^2} + C\int_0^t\ll(1 + \esp{X_s^p} + \esp{Y_s^p}\rr)ds \\
	+ C\ll(\int_0^t\ll(1+ \esp{X_s^{2p}} + \esp{Y_s^{2p}}\rr)ds\rr)^{1/2}\\
	\leq  \esp{X_0^p} + C\int_0^t\ll(1+\esp{X_s^{2p}} + \esp{Y_s^{2p}}\rr)ds\\
	\leq C_t\ll(1+\esp{X_0^{2p}} + \esp{Y_0^{2p}}\rr).
	\end{multline*}
	The result for the process~$Y$ follows from the same computation using~\eqref{eqy} instead of~\eqref{eqx}.
\end{proof}

{
	\begin{remark}
		In the proof of Lemma~\ref{aprioriestimate}, we use the stopping times $\tau_M$ to guarantee the function $u_{M,p}$ to be locally integrable. Otherwise we would not have the right to apply Gr\"onwall's lemma.
	\end{remark}
	
	The following lemma is similar to the previous one: it allows to control the moments of the process $(X^N,Y^N)$ instead of the moments of the limit process $(\bar X,\bar Y).$ An important point of the lemma is that, the desired estimates are uniform in~$N$ (provided that the initial conditions~$(X^N_0,Y^N_0)$ admit finite moments uniformly bounded in~$N$). This lemma is used in the proof of Theorem~\ref{mainresult} at Section~\ref{proofmainresult}. In the proof of Lemma~\ref{aprioriXNYN} below, we do not introduce the stopping times $\tau_M$ as in the proof of Lemma~\ref{aprioriestimate} for the sake of readability.
}

	\begin{lemma}\label{aprioriXNYN}
		Grant Assumption~\ref{hyp}$.(i)$ and Assumption~\ref{hyp}$.(iv).$ Let $X^N$ be the solution of the first equation of~\eqref{XNYN}. Assume that each $X^N_0$ has a finite second order moment. Then, for any $N\in\n^*$ and $t\geq 0,$
		$$\underset{s\leq t}{\sup}~\esp{\ll(X^N_s\rr)^2}\leq C_t\ll(1 +\esp{(X^N_0)^2}\rr),$$
		where $C_t>0$ does not depend on~$N.$
	\end{lemma}
	
	\begin{proof}
		By Ito's formula,
		\begin{align*}
		(X^N_t)^2 =& (X^N_0)^2 - 2\alpha_1\int_0^t (X^N_s)^2ds \\
		&+ \sum_{j=1}^N\int_{[0,t]\times\r_+\times\r} \ll[\ll(X^N_{s-} + \frac{u}{\sqrt{N}} - \frac{X^N_{s-}}{N}\rr)^2 - (X^N_{s-})^2\rr]\uno{z\leq f_1(X^N_{s-})}d\pi^{(1),j}(s,z,u)\\
		\leq& (X^N_0)^2 \\&+ \sum_{j=1}^N\int_{[0,t]\times\r_+\times\r} \ll[\ll(X^N_{s-} + \frac{u}{\sqrt{N}} - \frac{X^N_{s-}}{N}\rr)^2 - (X^N_{s-})^2\rr]\uno{z\leq f_1(X^N_{s-})}d\pi^{(1),j}(s,z,u)\\
		\leq& (X^N_0)^2 + \sum_{j=1}^N\int_{[0,t]\times\r_+\times\r} \ll[2\ll(\frac{u}{\sqrt{N}} - \frac{X^N_{s-}}{N}\rr)X^N_{s-} + \ll(\frac{u}{\sqrt{N}} - \frac{X^N_{s-}}{N}\rr)^2\rr]\\
		&\quad\quad\quad\quad\quad\quad\times\uno{z\leq f_1(X^N_{s-})}d\pi^{(1),j}(s,z,u).
		\end{align*}

		Taking expectation, we obtain
		\begin{align*}
		\esp{(X^N_t)^2}\leq &\esp{(X^N_0)^2} \\
		&+ N||f_1||_\infty\int_0^t\int_\r \esp{2\ll(\frac{u}{\sqrt{N}} - \frac{X^N_{s}}{N}\rr)X^N_{s} + \ll(\frac{u}{\sqrt{N}} - \frac{X^N_{s}}{N}\rr)^2}d\nu(u)ds.
		\end{align*}
		
		Now let us expand the term within the expectation
		$$2\ll(\frac{u}{\sqrt{N}} - \frac{X^N_{s}}{N}\rr)X^N_{s} + \ll(\frac{u}{\sqrt{N}} - \frac{X^N_{s}}{N}\rr)^2=2\frac{u}{\sqrt{N}}X^N_s - 2\frac{(X^N_s)^2}{N} + \frac{u^2}{N} + \frac{(X^N_s)^2}{N^2} - 2\frac{u X^N_s}{N\sqrt{N}}.$$
		
		Recalling that $\int ud\nu(u)=0$ and that $\sigma^2 := \int u^2d\nu(u),$ we have
		\begin{multline*}
		\int_\r \esp{2\ll(\frac{u}{\sqrt{N}} - \frac{X^N_{s}}{N}\rr)X^N_{s} + \ll(\frac{u}{\sqrt{N}} - \frac{X^N_{s}}{N}\rr)^2}d\nu(u) \\
		=  -\frac{2}{N}\esp{(X^N_s)^2} + \frac{\sigma^2}{N} + \frac{1}{N^2}\esp{(X^N_s)^2}\\
		\leq \frac{\sigma^2}{N} + \frac{1}{N^2}\esp{(X^N_s)^2}.
		\end{multline*}
		
		Consequently,
		$$\esp{(X^N_t)^2}\leq \esp{(X^N_0)^2} + \sigma^2||f_1||_\infty t + \frac1N\int_0^t\esp{(X^N_s)^2}ds.$$
		
		And Gr\"onwall's lemma implies that
		$$\esp{(X^N_t)^2}\leq \ll(\esp{(X^N_0)^2} + \sigma^2||f_1||_\infty t\rr)e^{t/N}\leq \ll(\esp{(X^N_0)^2} + \sigma^2||f_1||_\infty t\rr)e^{t},$$
		which proves the result.
	\end{proof}

%%%%%%%%%%%%%%%%%%%%%%%%%%%%%%%%%%%%%%%%%%%%%%%%%%%%%%%%%%%%%%%%%%%
%%                                                               %%
%% You may add acknowledgments (optional).                       %%
%%                                                               %%
%%%%%%%%%%%%%%%%%%%%%%%%%%%%%%%%%%%%%%%%%%%%%%%%%%%%%%%%%%%%%%%%%%%

%\ACKNO{We are grateful to the authors of the ECP/EJP class.}

%%%%%%%%%%%%%%%%%%%%%%%%%%%%%%%%%%%%%%%%%%%%%%%%%%%%%%%%%%%%%%%%%%%
%%                                                               %%
%% You have reached the end of your document.                    %%
%%                                                               %%
%%%%%%%%%%%%%%%%%%%%%%%%%%%%%%%%%%%%%%%%%%%%%%%%%%%%%%%%%%%%%%%%%%%
\nocite{*}
\end{document}